\documentclass[a4paper,12pt]{article}
\usepackage{amsmath,latexsym,amssymb,enumerate,newcom}
\usepackage[mathscr]{eucal}
\usepackage[dvips]{graphicx}
\begin{document}
\title{\sf{A convergence rate of
           the discrepancy for Allen-Cahn type equations}}
\author {\sf{Kazuhiro HORIHATA}}
\maketitle
\begin{abstract}
{\bf{abstract}.} \enspace
This paper presents a new partial differential equation to build Brakke's motion, 
which is a weak notion of mean curvature flow.
We call the equation a modified Allen-Cahn equation abbreviated to \textit{MAC}.
After that we introduce its benefit:
An improved estimate on the discrepant energy that means the equipartition of the first energy.
This equation is obtained by rendering the parameter in Allen-Cahn equation 
time-dependent and adding a term like Z-transform.
Furthermore, we mention the existence of a Brakke's motion through our \textit{MAC}.
\end{abstract}
%
%
\numberwithin{equation}{section}
\section{\enspace Introduction}
\par
If a smooth family of embedded hypersurfaces $\Gamma_t\,\subset\,\mathbb{R}^{N+1}$ $(t\in I)$ follows%
\begin{equation}\partial_t{x}\;=\;{\mathbf{H}}(x)\label{EQ:MCF}\end{equation}
for \begin{math}{x}\,=\,{x}(t)\,\in\,\Gamma_t\end{math} and $t\in I$,
$\{\Gamma_t\}$ is called \lq\lq mean curvature motion.\rq\rq
Here \begin{math}I\,\subset\,\mathbb{R}_+\end{math} is a right semi-open interval,
$\partial_t{x}$ the normal velocity at $x$ and 
${\mathbf{H}}(x)$ the mean curvature vector at $x$.
Typical examples are
\begin{Exa}\label{EXA:ShrinkingSpheresandCylidenrs}
If \begin{math}\Gamma_t\,=\,\partial B_{r(t)}^{N+1}\,\subset\,\mathbb{R}^{N+1},\end{math}
then \eqref{EQ:MCF} reduces to an ordinary differential equation for the radius{\rm{:}}
\begin{equation*}\dot{r}\,=\,-\frac Nr.\end{equation*}
The solution with $r(0)=R$ is \begin{math}r(t)\,=\,\sqrt{R^2\,-\,2Nt}\end{math}
$(t\in (-\infty,R^2/(2N)))$. 
Similarly we have the shrinking cylinders 
\begin{math}\Gamma_t\,=\,\mathbb{R}^j\times\partial B_{r(t)}^{N+1-j}\,\subset\,\mathbb{R}^{N+1}\end{math}
$(j=1,2,\ldots,N)$ with
\begin{math}r(t)\,=\,\sqrt{R^2\,-\,2(N-j)t}\end{math} $(t\in (-\infty,R^2/(2(N-j))))$.
\end{Exa}
\begin{Exa}\label{EXA:GrimReaper}
For $N=1$, an explicit solution given by $\Gamma_t\,=\,\graph u_t$
with \begin{math}u_t(p)\,=\, t-\log\cos p\end{math} $(p\in (-\pi/2,\pi/2))$
is called \lq\lq the grim reaper solution.\rq\rq
\end{Exa}
See Ecker\cite[Chap.~2]{ecker}.
\par
Instead of viewing the mean curvature flow as an evolution equation for the hypersurfaces $\Gamma_t$,
we can also regard it as an evolution equation for a smooth family of embedding
\begin{math}X=X(\cdot,t):\mathbb{R}^N\times I\,\to\,\mathbb{R}^{N+1}\end{math}
with \begin{math}\Gamma_t\,=\,X(M,t)\end{math} for some \begin{math}M\,\subset\,\mathbb{R}^N.\end{math}
Setting $x=X(p,t)$, equation \eqref{EQ:MCF} takes the form
\begin{equation}\frac{\partial X}{\partial t}(p,t)\;=\;\triangle_{\Gamma_t} X(p,t).
\label{EQ:MCF-X}\end{equation}
\par If $\{\Gamma_t\}$ $(t\in I)$ follows its smooth mean curvature, 
it has the integral form: $\{\Gamma_t\}$ yields
\begin{equation}
\frac d{dt} \lint_{\Gamma_t}\phi\, d\mathcal{H}^{N}
\;=\;-\lint_{\Gamma_t}\phi|{\mathbf{H}}|^2\, d\mathcal{H}^{N}
\,+\,\lint_{\Gamma_t} D\phi\cdot{\mathbf{H}}\, d\mathcal{H}^{N}
\label{EQ:IntOfMCM}\end{equation}
for any \begin{math}\phi\,\in\,C_c^2(\mathbb{R}^{N+1}\times I).\end{math}
In 1978, Brakke\cite{brakke} gave a weak notion for \eqref{EQ:IntOfMCM},
furthermore nearly two decades years later Ilmanen\cite{ilmanen} has improved it as follows: 
\begin{Defn}\label{DEF:BrakkeMotion}
Let $\mu$ be a Radon measure and $\phi$ be of $C_c^2 (\mathbb{R}^{N+1}\times\mathbb{R}_+)$.
Then $\beta (\mu, \phi)$ be given by 
\begin{equation}+\infty\;\mathrm{if}\;
\mu\lfloor_{\{\phi > 0\}}\;\not\in\,\mathscr{M}_N\;\mathrm{or}\;
\vert\delta\mu\vert\lfloor_{\{\phi > 0\}}\, \not\ll\,\mu\lfloor_{\{\phi > 0\}}\;\mathrm{or}\;
\int\phi|\mathbf{H}|^2\, d\mu\,=\,+\infty
\notag\end{equation}
where $\mathscr{M}_N$ is the set of all $N$-rectifiable Radon measures.
Otherwise
\begin{equation}
\beta (\mu,\phi)\;=\;\int-\phi|{\mathbf{H}}|^2
\,+\, D\phi\cdot\bigl((T_x\mu)^\bot\cdot{\mathbf{H}}\bigr)\, d\mu.
\notag\end{equation}\end{Defn}
So Brakke's motion $\{\mu_t\}$ $(t\in I)$ is indicated by
\begin{equation}\overline{D}_t\mu_t (\phi)\;\le\;\beta(\mu_t,\phi)\label{INEQ:Brakke}\end{equation}
for all $t\ge 0$. Here
\begin{equation}\overline{D}_t\mu_t\;=\;\limsup_{s\to t}\frac{\mu_s\,-\,\mu_t}{s-t},\notag\end{equation}
$T_x\mu$ is the approximate tangent space and $(T_x\mu)^\perp$ is the perpendicular space of it and
$\mathbf{H}$ is generalized mean curvature. 
See Simon\cite[\S 16.5]{simon}.
\par Another approach to the mean curvature flow is to consider a signed distance function given by
\allowdisplaybreaks\begin{align}
& d(x,t):\mathbb{R}^{N+1}\times I\,\to\,\mathbb{R},\notag\\
& M_0\;=\;\{x\,\in\,\mathbb{R}^{N+1}\,;\, d(x,0)\,>\,0\},\; 
M_t\;=\;\{x\,\in\,\mathbb{R}^{N+1}\,;\, d(x,t)\,>\,0\},\notag\\
& \Gamma_0\;=\;\{x\,\in\,\mathbb{R}^{N+1}\,;\, d(x,0)\,=\,0\},\notag\\
& \Gamma_t(\,=\,\partial M_t)\,=\,\{x\,\in\,\mathbb{R}^{N+1}\,;\, d(x,t)\,=\,0\},\;\notag\\
& |D d|\;=\;1\quad\mathrm{near}\;\Gamma_t,\;
d_t\;=\;\triangle d\quad\mathrm{on}\;\Gamma_t.\label{EQ:Dis}\end{align}
Note that the normal velocity of the motion of the hypersurfaces $\Gamma_t$ at $x$ is $-d_t$ and 
the mean curvature $H$ $-\triangle d$.
For the above \eqref{EQ:Dis}, we know
\begin{Prop}\label{PROP:DistanceFctOfMCF}
\begin{enumerate}
\item[\rm{(1)}] \begin{math}\Gamma_0\,\in\, C^{k+\alpha}\, (k\ge 2, 0<\alpha<1)\end{math} is 
the boundary of a bounded open set, 
then there exists a unique mean curvature flow $\Gamma_t$ in $I$ starting from $\Gamma_0$\rm{;} 
the signed distance function $d(x,t)$ belongs to $C^{k+\alpha,(k+\alpha)/2}$
in a small neighborhood of \begin{math}\bigcup_{t\in I} (\Gamma_t\times\{t\}).\end{math}
\item[\rm{(2)}] If $N=1$ or $\Gamma_0$ is convex, then the solution can be extended up to a time at 
which $\Gamma_t$ shrinks to a single point.
\end{enumerate}\end{Prop}
\par (1) was the result by Hamilton \cite{hamilton} 
and also by Evans-Spruck \cite{evans-spruck-1}. 
While (2) was proved by Gage-Hamilton \cite{gage-hamilton} if $N=1$ and $\Gamma_0$ is convex,
by Grayson \cite{grayson} in $N=1$
and by Huisken \cite{huisken} provided that $N\ge 1$ and $\Gamma_t$ is convex.
\par Viscosity solution, which is another weak global solution, was established 
by Evans-Spruck \cite{evans-spruck-2},\cite{evans-spruck-3}
and Chen-Giga-Goto \cite{chen-giga-goto}.
\par Ilmanen \cite{ilmanen} constructed a Brakke's motion from Allen-Cahn equation described as
\begin{equation}
\frac{\partial u_\epsilon}{\partial t}\,-\,\triangle u_\epsilon\,+\,
\frac 2{\epsilon^2} (u_\epsilon^2\,-\,1)u_\epsilon\;=\;0.
\label{EQ:Allen-Cahn}\end{equation}
The effect of \begin{math}(u_\epsilon^2-1)u_\epsilon/\epsilon^2\end{math} is to force $u_\epsilon$
to approximate a step function.
\par Changing an unknown function \begin{math}u_\epsilon\;=\;q(z_\epsilon/\epsilon)\end{math}
with \begin{math}q\;=\;\tanh,\end{math} \eqref{EQ:Allen-Cahn} becomes
\begin{equation}
\frac{\partial z_\epsilon}{\partial t}\,-\,\triangle z_\epsilon\,+\,
\frac {2q}{\epsilon} (|D z_\epsilon|^2\,-\,1)\;=\;0.
\label{EQ:ApproDis}\end{equation}
\eqref{EQ:ApproDis} can be regarded as an approximation of \eqref{EQ:Dis}
because it becomes the linear heat flow if \begin{math}|D z_\epsilon|=1\end{math}
or \begin{math} z_\epsilon=0.\end{math}
To briefly explain his paper, we prepare a few notions:
the energy equality called $\epsilon$-Brakke formula includes
\begin{equation}
\frac d{dt}\lint\phi\, d\mu_t^\epsilon
\;=\;-\epsilon\lint\phi\left|\frac{\partial u_\epsilon}{\partial t}\right|^2 \, dx
\,-\,\lint S:D^2\phi\, dV_t^\epsilon\,-\,\lint\nu\otimes\nu:D^2\phi\, d\xi_t^\epsilon,
\label{EQ:EBrakke}\end{equation}
where
\allowdisplaybreaks\begin{align}
& d\mu_t^\epsilon\,=\,\left(\frac\epsilon 2|Du_\epsilon|^2\,+\,\frac 1{2\epsilon}(u_\epsilon^2-1)^2\right)
\, dx, \;
d\xi_t^\epsilon\,=\,\left(\frac 1{2\epsilon}(u_\epsilon^2-1)^2\,-\,\frac\epsilon 2|Du_\epsilon|^2\right)
\, dx,\notag\\
& S\,=\,(S_{i,j}), \; S_{i,j}\,=\,\delta_{i,j}\,-\,\nu_i\nu_j,\; \nu_i
\,=\,\frac{D_i u_\epsilon}{|Du_\epsilon|},\notag\\
& V_t^\epsilon(\varPsi)\;=\;\lint \varPsi(x,Du_\epsilon^\perp)\, d\mu_t^\epsilon
\quad (\varPsi\,\in\,C_c^0(\mathbb{R}^{N+1}\times G_{N}(\mathbb{R}^{N+1}))).\end{align}
\par His paper \cite{ilmanen} is heuristically composed of three parts:
The first demonstrates a few energy equalities or inequalities, 
including the approximation flow of \eqref{EQ:EBrakke} 
and by using it, the $\epsilon$-Brakke's version of Huisken's monotonicity formula.
The second discusses both of a few local properties, e.g., clearing-out lemma 
and a global property, e.g., the discrepancy theorem.
The third constructs Brakke's motion, combining the second step 
with rectifiability theorem and lower-semicontinuity by Allard \cite{allard} and Brakke \cite{brakke}.
\par This paper, proposes a new P.D.E to build mean curvature flow:
\begin{equation}
\frac{\partial v_\epsilon}{\partial t}\,-\,\triangle v_\epsilon\,+\,
\frac 2{\epsilon^{2(1-\kappa)}} (v_\epsilon^2\,-\,1)v_\epsilon
\,-\,\frac{\dot{\kappa}}{K(\epsilon)}(v_\epsilon^2\,-\,1)\log\frac{1+v_\epsilon}{1-v_\epsilon}
\;=\;0\label{EQ:Mac}\end{equation} and states a benefit of it.
Here \begin{math}\kappa(t)\,=\,\arctan (t)/\pi\end{math} and
\begin{math}K(\epsilon)\,=\,2/|\log\epsilon|.\end{math}
We say a new P.D.E \eqref{EQ:Mac} to be a modified Allen-Chan equation shorten it to \lq\lq\textit{MAC}:\rq\rq
By transformation of \begin{math}v_\epsilon (x,t)\;=\; q(r_\epsilon(x,t)/\epsilon^{1-\kappa}),\end{math}
we see that $r_\epsilon$ must satisfy
\begin{equation}
\frac{\partial r_\epsilon}{\partial t}\,-\,\triangle r_\epsilon
\,+\,\frac{q}{\epsilon^{1-\kappa}}(|Dr_\epsilon|^2\,-\,1)\;=\;0.
\label{EQ:NewDis}\end{equation}
Since we have
\begin{equation}
\frac{\partial v_\epsilon}{\partial t}\;=\;\frac{\dot{q}}{\epsilon^{1-\kappa}}
\left(\frac{\partial r_\epsilon}{\partial t}\,-\,|\log\epsilon|\dot{\kappa}r_\epsilon\right),
\notag\end{equation}
the third term in \eqref{EQ:Mac}
is necessary for \eqref{EQ:NewDis} to hold.
\par The rest of our paper has three folds:
The first commences with studying a few properties of our nonlinear term because it is probably new. 
The second, the main topic of this paper shows a classical existence of the solution,
an energy inequality and the discrepancy theorem that means the equipartition of 
\begin{equation}
\liint\, dtd\xi_t^\epsilon \;=\;\llint
\,\left(\frac{(v_\epsilon^2-1)^2}{2\epsilon^{1-\kappa}}\,-\,\frac{\epsilon^{1-\kappa}}2|Dv_\epsilon|^2\right)
\, dtdx\;=\;O(\frac 1{|\log\epsilon|})\quad (\epsilon\searrow 0).                                              
\label{EQ:Discrepancy}\end{equation}
Here notice that benefit of our equation is to be able to compute a convergence rate
of the discrepant energy, which means the second term of \eqref{EQ:Discrepancy}, 
all discussions here are elementary though.
The second introduces a comparison theorem, which plays a crucial role in proving clearing-out lemma 
in the construction of Brakke's motion via \textit{MAC}. 
We finally give a remark that our $\epsilon$-Brakke's formula so similar to \eqref{EQ:EBrakke} 
does converge to a Brakke's motion \eqref{INEQ:Brakke}.
Its proof is as same as the one by Ilmanen \cite{ilmanen}.

%
%
\section{A new P.D.E}\label{newPDE}
This section presents a new partial differential equation of Allen-Cahn type.
\subsection{Nonlinear term}
Set $\epsilon$ be any positive number{\rm{;}} We choose $\kappa$ and $K$ as
\allowdisplaybreaks\begin{align}
&\kappa (t)\;=\;\frac 1\pi\arctan(t)\; (t\,\in\,\mathbb{R}),\;
K(\epsilon)\,=\,2/\vert\log\epsilon\vert.
\notag\end{align}
When we define a few functions $f$, $g$ and $\phi_\epsilon$  and their integration by
\begin{equation}
f(x)\;=\;\bigl(x^2\,-\,1\bigr)x,\;                                        
g(x)\;=\;(x^2\,-\,1)\log\Bigl\vert\dfrac{1+x}{1-x}\Bigr\vert,\;
\phi_\epsilon(x)\;=\;\dfrac 2{\epsilon^{2(1-\kappa)}}f(x)\,-\,\dfrac{\dot{\kappa}}{K(\epsilon)} g(x).
\notag\end{equation}
\allowdisplaybreaks\begin{align}
& F(u)\,=\;\frac 12(u^2-1)^2,\notag\\
& G(u)\,=\,\frac 13(u-2)(u+1)^2\log\vert 1+u\vert\,-\,\frac 13(u+2)(u-1)^2\log\vert 1-u\vert
\,+\,\frac{u^2}3\,+\,\frac 43\log 2\,-\,\frac 13,\notag\\            
& W_\epsilon(u)\,=\,\frac {F(u)}{2\epsilon^{1-\kappa}}\,-\,\frac{\dot\kappa}{K(\epsilon)} G(u).
\notag\end{align}
The nonlinear term $\phi_\epsilon$ has the following properties:
there exist two positive functions $\alpha_\epsilon (t)$ and $\beta_\epsilon (t)$ 
\begin{math}(\alpha_\epsilon<1<\beta_\epsilon)\end{math} with
\allowdisplaybreaks\begin{align}
& \phi_\epsilon(\pm\beta_\epsilon)\,=\,\phi_\epsilon(\pm 1)\,=\,\phi_\epsilon(\pm\alpha_\epsilon)
\,=\,\phi_\epsilon(0)\,=\, 0,\notag\\ 
& \phi_\epsilon\,>\,0\quad\mathrm{on}\quad (-\beta_\epsilon,-1)\cup(-\alpha_\epsilon,0) 
\cup (\alpha_\epsilon,1)\cup (\beta_\epsilon,\infty),
\notag\\
&\phi_\epsilon \, < \, 0 \quad \mathrm{on}\quad (-\infty,-\beta_\epsilon)\cup (-1,-\alpha_\epsilon)
\cup (0,\alpha_\epsilon)\cup (1,\beta_\epsilon),
\notag\\
& \phi_\epsilon^\prime(\pm\beta_\epsilon)\,>\,0,\;\phi_\epsilon^\prime(\pm\alpha_\epsilon)\,>\, 0, 
\phi_\epsilon^\prime (\pm 1)\,<\, 0,\; \phi_\epsilon^\prime (0)\,<\,0,
\notag\\
&
g_\epsilon \,>\,0\quad\mathrm{on}\quad (-1,0)\cup (1,\infty),\; g_\epsilon \,<\,0
\quad\mathrm{on}\quad (-\infty,-1)\cup (0,1),\notag\\
&
g_\epsilon (0)\,=\,g_\epsilon (\pm 1)\,=\,0,\;g_\epsilon^\prime (\pm 1)\,=\,+\infty,\;g_\epsilon^\prime (0)\,<\,0.
\label{COND:Eq}\end{align}
\begin{Rem}
When we put
\allowdisplaybreaks\begin{align}                                                                               
& g_{\epsilon}\,=\,\frac{\vert\log\epsilon\vert}2\dot{\kappa}\epsilon^{2(1-\kappa)},\notag\\
& 1-\delta_{\epsilon}\,=\,\alpha_\epsilon\,=\,g_{\epsilon}\log\frac{1+\alpha_\epsilon}{1-\alpha_\epsilon},\;   
1+\eta_{\epsilon}\,=\,\beta_\epsilon\,=\,g_{\epsilon}\log\frac{\beta_\epsilon+1}{\beta_\epsilon-1},
\notag\end{align}
then $\delta_{\epsilon}$ and $\eta_{\epsilon}$ satisfy
\allowdisplaybreaks\begin{align}
& 2\exp\bigl(-\frac 1{|\log\epsilon|\dot\kappa\epsilon^{2(1-\kappa)}}\bigr)\,<\,\delta_{\epsilon}
\,<\,4\exp\bigl(-\frac 1{|\log\epsilon|\dot\kappa\epsilon^{2(1-\kappa)}}\bigr),\notag\\
& 2\exp\bigl(-\frac 4{|\log\epsilon|\dot\kappa\epsilon^{2(1-\kappa)}}\bigr)\,<\,\eta_{\epsilon}
\,<\,4\exp\bigl(-\frac 2{|\log\epsilon|\dot\kappa\epsilon^{2(1-\kappa)}}\bigr).
\notag\end{align}\end{Rem}
\begin{Rem}
Set \begin{math}t_\epsilon\,=\,|\log\epsilon|/\pi.\end{math}
The points $\alpha_\epsilon(t)$ and $\beta_\epsilon(t)$ $(t>0)$ behave the following\rm{:}
$\alpha_\epsilon(t)$ decreases in $0\le t<t_\epsilon$ and increases in $t_\epsilon$ $<$ $t$,
while $\beta_\epsilon(t)$ increases in \begin{math}0\le t<t_\epsilon\end{math}
and decreases in $t_\epsilon<t$.\end{Rem}
\begin{Rem}\label{REM:2-1-1}
From the symmetricity of $\phi_\epsilon$, we say
\allowdisplaybreaks\begin{align}
&
\int_{-\gamma_\epsilon}^0\phi_\epsilon(t)\, dt\,+\,\int_{0}^{\gamma_\epsilon}\phi_\epsilon(t)\, dt\,=\, 0
\notag 
\label{EQ:2-1-1}\end{align}\end{Rem}
where $\gamma_\epsilon$ is $\alpha_\epsilon$ or $1$ or $\beta_\epsilon$.
\subsection{A new P.D.E}
Now we are in the position to introduce
a new partial differential equation 
and a benefit of it:
\begin{equation}
\renewcommand{\arraystretch}{1.2}
\left\{
\begin{array}{ll}
\dfrac{\partial v_\epsilon}{\partial t}\,-\,\triangle v_\epsilon\,+\,\phi_\epsilon(v_\epsilon)\;=\;0 
&\quad\mathrm{in}\quad (0,t_\epsilon]\times\mathbb{R}^{N+1}, \notag\\
v_\epsilon\,=\, g &\quad\mathrm{in}\quad\{0\}\times\mathbb{R}^{N+1}\notag
\end{array}\right.\tag{MAC}\label{EQ:P}\end{equation}
where $g$ is a bounded and $C^{2+\alpha}$-function defined on $\mathbb{R}^{N+1}$ $(0<\alpha<1)$.
We need either of the hypotheses below for each of our theorems to hold:
\allowdisplaybreaks\begin{align}
&\vert g\vert\,\le\,\alpha_\epsilon(t_\epsilon)\quad
\mathrm{and}\quad \vert g\vert\,=\,\alpha_\epsilon(t_\epsilon)\; (|x|>K_0)\;
\text{for some positive}\; K_0,\label{INEQ:I1}
\tag{$\mathrm{I}_1$}\\
& 
Dg\,\in\, L^2, W_\epsilon(g)\,\in\, L,\label{INEQ:I2}
\tag{$\mathrm{I}_2$}\\
& \frac 1{2\epsilon}\bigl(g^2\,-\,1 \bigr)^2\,-\,\frac {\epsilon}2\vert Dg\vert^2\;\ge\; 0\label{INEQ:I3}
\tag{$\mathrm{I}_3$}.\end{align} 
The first energy and the discrepant energy density of \eqref{EQ:P} is respectively given by
\allowdisplaybreaks\begin{align}
& d\mu_t^\epsilon\,=\,
\Bigl(\frac{\epsilon^{1-\kappa}}2\vert Dv_\epsilon\vert^2\,+\,\frac{F(v_\epsilon)}{2\epsilon^{1-\kappa}}
\,-\,\frac{\epsilon^{1-\kappa}\dot\kappa}{K(\epsilon)}G(v_\epsilon)\Bigr)\,dx,
\label{DEF:FE}\\
& d\xi_t^\epsilon\,=\,\Bigl(\frac{F(v_\epsilon)}{2\epsilon^{1-\kappa}}
\,-\,\frac{\epsilon^{1-\kappa}}2\vert Dv_\epsilon\vert^2\Bigr)\, dx,\label{DEF:Dis}\end{align}
by transformation of
\begin{equation}v_\epsilon\,=\,q\Bigl(\frac{r_\epsilon}{\epsilon^{1-\kappa}}\Bigr)\notag\end{equation}
\eqref{DEF:Dis} is expressed by
\begin{equation} 
d\xi_t^\epsilon\,=\,\frac {F(v_\epsilon)}{2\epsilon^{1-\kappa}}\,(1\,-\,\vert Dr_\epsilon\vert^2)\, dx. 
\label{EQ:Dis-2}\end{equation}
\setcounter{chapternumber}{3}\setcounter{equation}{0}
\renewcommand{\theequation}%
             {\thechapternumber.\arabic{equation}}
\section{\enspace Existence, Energy equality and Discrepancy}
This section describes our main theorems. First
\begin{Thm}{\rm{(Existence).}}\label{THM:Exist}
Let us assume $(\mathrm{I}_1)$ for the initial value $g$\rm{;}
MAC has a classical solution
with \begin{math}|v_\epsilon|\,\le\,\alpha_\epsilon (t_\epsilon).\end{math}
\end{Thm}
In the following theorem, let $v_\epsilon$ be a solution to \eqref{EQ:P}:
\begin{Thm}{\rm{(Energy Conservative Law).}}\label{THM:EEq} 
Suppose that the initial condition $g$ satisfies \eqref{INEQ:I1} and \eqref{INEQ:I2}.
Then we claim
\allowdisplaybreaks\begin{align}
& \int_0^t\epsilon^{1-\kappa}\,\lint_{\mathbb{R}^{N+1}}
\left|\frac{\partial v_\epsilon}{\partial t}\right|^2\, dx \,+\,\lint_{\mathbb{R}^{N+1}}\, d\mu_t^{\epsilon} 
\label{EQ:EI}\\
&\,+\,\int_0^t\frac{\dot\kappa\, ds}{K(\epsilon)}\,\lint_{\mathbb{R}^{N+1}}\, d\xi_s^\epsilon
\,+\,\int_0^t \frac{\ddot\kappa\,+\,\dot\kappa^2|\log\epsilon|}{K(\epsilon)\dot\kappa}
\epsilon^{1-\kappa} \, ds \lint_{\mathbb{R}^{N+1}} G(v_\epsilon) \, dx
\;=\;\lint_{\mathbb{R}^{N+1}}\, d\mu_0^{\epsilon}\notag\end{align}
for all \begin{math}t\,\in\,[0,t_\epsilon].\end{math}\end{Thm}
From Theorem \ref{THM:EEq}, we directly obtain
\begin{Thm}{\rm{(Discrepancy theorem)}}\label{THM:Main}
Under the hypothesis \eqref{INEQ:I1},\eqref{INEQ:I2} and \eqref{INEQ:I3}, the following holds
\begin{equation}
0\;\le\;\int_0^{t_\epsilon}{\dot\kappa (t)}\, dt\lint_{\mathbb{R}^{N+1}}\, d\xi_t^{\epsilon}
\;=\; O\left(\frac 1{|\log{\epsilon}|}\right)\quad (\epsilon\searrow 0).
\end{equation}\end{Thm}
\vskip 9pt\noindent{\underbar{{Proof of Theorem \ref{THM:Exist}}}.}
\rm\enspace\vskip 6pt
First of all, note that 
a maximal principle leads us to \begin{math}|v_\epsilon|\,\le\,\alpha_\epsilon (t_\epsilon)\end{math}
if the solutions are smooth enough.
Next we shall prove a classical existence.
Set any number \begin{math}t_1\,=\,\theta_0\epsilon^2\end{math}
with \begin{math}1/\theta_0\;=\;\sqrt{2}\sqrt{1+\sqrt{2}}\,+\,\sqrt{2(1+\sqrt{2})}/e\end{math}
and $K$ a positive integer greater than $K_0$,
and let $h_K$ be the solution to
\begin{equation}\left\{\begin{array}{lll}
& \dfrac{\partial h_K}{\partial t}\,-\,\triangle h_K\;=\; 0
&\mathrm{in}\;B_K^{N+1}\times (0,t_\epsilon],\\[2.5mm]
& \dfrac{\partial h_K}{\partial \nu}\;=\; 0
&\mathrm{on}\;\partial B_K^{N+1}\times (0,t_\epsilon],\\[2.5mm]
& h_K(0)\;=\; g&\mathrm{at}\;B_K^{N+1}.
\end{array}\right.\end{equation}
Then after setting \begin{math}w_{\epsilon,K}^{(0)}\;=\; 0\end{math} and 
\begin{math}v_{\epsilon,K}^{(0)}\,=\, h_K,\end{math}
we inductively define $w_{\epsilon,K}^{(n)}$ and $v_{\epsilon,K}^{(n)}$ $(n\in\mathbb{N})$ by
\allowdisplaybreaks\begin{align}
& w_{\epsilon,K}^{(n)}(x,t)\;=\;\,-\,\int_{0}^t\dfrac{2 ds}{\epsilon^{2(1-\kappa)}}
\,\lint_{B_K} G(t-s,x-y)\dot\chi\bigl(({v_{\epsilon,K}^{(n-1)}}^2\,-\,1)^2\bigr)
\bigl({v_{\epsilon,K}^{(n-1)}}^2\,-\,1\bigr) v_{\epsilon,K}^{(n-1)}(y,s)\, dy\notag\\ 
&\,+\,\lint_{0}^t\dfrac {\dot{\kappa}}{K(\epsilon)} \, 
\lint_{B_K}\, G(t-s,x-y)\dot\chi\bigl(({v_{\epsilon,K}^{(n-1)}}^2\,-\,1)^2\bigr)
({v_{\epsilon,K}^{(n-1)}}^2\,-\,1)
\log\biggl\vert\dfrac{1+v_{\epsilon,K}^{(n-1)}}{1-v_{\epsilon,K}^{(n-1)}}\biggr\vert (y,s)\, dy
\label{EQ:Duhamel1}\\
\intertext{and}
& v_{\epsilon,K}^{(n)}\;=\;h_K\,+\, w_{\epsilon,K}^{(n)}.
\label{EQ:Duhamel2}\end{align}
Here $G$ is the Green function for \begin{math}\partial/\partial t\,-\,\triangle\end{math}
in \begin{math}B_K^{N+1}\times(0,\infty)\end{math} with Neumann zero boundary,
moreover a smooth function $\chi$
is chosen by \begin{math}-3/2\le\chi\le 3/2\end{math} and
\begin{math}\dot\chi(\tau)\,=\,1\; (|\tau|\le 1),\;\,=\,0\; (|\tau|\ge 2).\end{math} 
We readily check that $|w_{\epsilon,K}^{(n)}|$ $\le$ $1$ and
\begin{math}w_{\epsilon,K}^{(n)}\,\in\, C^{2+\alpha,(2+\alpha)/2}(B_K\times [0,t_1])\end{math}
\begin{math}(0<\alpha<1)\end{math}
from Lady\v{z}henskaya,~O.~A., Solonnikov,~V.~A.,Ural'ceva,~N.~N
\cite[p.320, Theorem 5.2]{ladyzhenskaya-solonnikov-uralceva} because
\begin{equation}
(x^2\,-\,1)\log\left\vert\dfrac {1+x}{1-x}\right\vert
\notag\end{equation} is H\"{o}lder continuous.
Ascoli-Arzella's theorem tells us that by 
\begin{math}\dot\chi(\tau)\;=\;0\;\mathrm{if}\; \tau\ge 2,\end{math}
a subsequence of $\{w_{\epsilon,K}^{(n)}\}$ and $\{v_{\epsilon,K}^{(n)}\}$ $(n\in\mathbb{N})$ 
still denoted by them, respectively converges uniformly to 
\begin{math}C^{2+\alpha,(2+\alpha)/2}\end{math}-function $v_{\epsilon,K}$ and $w_{\epsilon,K}$ 
in $\overline{B}_K$ $\times$ $[0,t_1]$. 
Furthermore they satisfy \eqref{EQ:Duhamel1} by substituting $v_{\epsilon,K}$ and $w_{\epsilon,K}$ 
for $v_{\epsilon,K}^{(k)}$ and $w_{\epsilon,K}^{(k)}$ $(k=n-1,n)$.
\par By repeating the argument above for time-variable from $[0,t_1]$ to 
$[(j-1)t_1,jt_1]$ $(j=1,2,\cdots,[t_\epsilon/t_1]+1)$,
we conclude our assertion.
\par Next
we extend $v_{\epsilon,K}$ to the one in \begin{math}\mathbb{R}^{N+1}\times [0,t_\epsilon]\end{math}: 
Choose
\begin{equation}\eta_K(x)\;=\;\left\{\begin{array}{lll}
1& (|x|<K)\\0& (|x|>K+1)
\end{array}\right.
\quad |D\eta_K|\,\le\,2,\;|D^2\eta_K|\,\le\,4;\end{equation}
Then set an extension of $v_{\epsilon,K}$ by
\begin{equation}v_{\epsilon,K}\;=\;\left\{\begin{array}{ll}
v_{\epsilon,K}((|x|,x/|x|),t)&\; (|x|<K),\\
v_{\epsilon,K}((2K-|x|,x/|x|),t)\eta_K(x)&\; (K<|x|<2K),\\
0 & (2K<|x|).
\end{array}\right.\end{equation}
Furthermore we need the following:
\allowdisplaybreaks\begin{align}
& d\mu_t^{\epsilon,K}\,=\,\left(\frac{\epsilon^{1-\kappa}} 2|Dv_{\epsilon,K}|^2
\,+\,\frac 1{2\epsilon^{1-\kappa}}(v_{\epsilon,K}^2-1)^2\right)
\, dx,\notag\\
& d\xi_t^{\epsilon,K}\,=\,\left(\frac 1{2\epsilon^{1-\kappa}}(v_{\epsilon,K}^2-1)^2
\,-\,\frac{\epsilon^{1-\kappa}} 2|Dv_{\epsilon,K}|^2\right)
\, dx.\notag\end{align}
Note the extensive function is of $C_c^{2+\alpha,(2+\alpha)/2}$.
\par By a diagonal argument, we find that there exists a subsequence $\{K(j)\}$ $(j\in\mathbb{N})$
of $\{K\}$ such that $v_{\epsilon,K(j)}$ locally converges to some $v_\epsilon$ 
in the sense of $C^{2+\alpha,(2+\alpha)/2} (\mathbb{R}^{N+1}\times[0,t_{\epsilon}])$.
Since $v_{\epsilon,K(j)}$ satisfies \eqref{EQ:Mac} 
in \begin{math}B_{K(j)}\times [0,t_\epsilon],\end{math}
we can deduce that $v_\epsilon$ does \eqref{EQ:Mac} in \begin{math}B_{L}\times [0,t_\epsilon]\end{math}
for any positive integer \begin{math}L\end{math} greater than $K_0$.
\vskip 9pt
\noindent{\underbar{{Proof of Theorem \ref{THM:EEq}}.}
\rm\enspace
\vskip 6pt
First recall that $v_{\epsilon,K}$ implies 
\begin{equation}\left\{\begin{array}{lll}
& \dfrac{\partial v_{\epsilon,K}}{\partial t}\,-\,\triangle v_{\epsilon,K}
\,+\,\phi_\epsilon (v_{\epsilon,K})\;=\; 0
&\mathrm{in}\;B_K^{N+1}\times (0,t_\epsilon],\\[2.5mm]
& \dfrac{\partial v_{\epsilon,K}}{\partial \nu}\;=\; 0
&\mathrm{on}\;\partial B_K^{N+1}\times (0,t_\epsilon],\\[2.5mm]
& v_{\epsilon,K}(0)\;=\; g&\mathrm{at}\;B_K^{N+1},
\end{array}\right.\label{EQ:Mac_K}\end{equation}
multiply \begin{math}\partial v_{\epsilon,K}/\partial t\end{math} by \eqref{EQ:Mac_K} and
integrate it on \begin{math}B_K\times (0,t)\end{math} for any positive time
\begin{math}t\,\le\,t_\epsilon\end{math} to verify
\allowdisplaybreaks\begin{align}
&\lint_0^t\epsilon^{1-\kappa}\,ds \lint_{B_K}\left|\frac{\partial v_{\epsilon,K}}{\partial s}\right|^2\, dx
\,+\,\lint_{B_K}\,d\mu_t^{\epsilon,K}
\,+\,\lint_0^t \frac{\dot\kappa\, ds}{K(\epsilon)}\,\lint_{B_K}\, d\xi_t^{\epsilon,K}
\label{EQ:EnergyEQ}\\                
& \,+\,\int_0^t \frac{\ddot\kappa\,+\,\dot\kappa^2|\log\epsilon|}{K(\epsilon)\dot\kappa}
\epsilon^{1-\kappa}\, ds\lint_{B_K} G(v_{\epsilon,K})\, dx\;=\;\lint_{B_K}\, d\mu_0^{\epsilon}.
\notag\end{align}
\par Thus, on account of $d\xi_t^{\epsilon,K}$ $\ge$ $0$ and 
${\ddot\kappa}/{\dot\kappa}\,+\,\dot\kappa|\log\epsilon|\,\ge\, 0$
in $[0,t_\epsilon]$,
by monotone convergence theorem, 
we pass to the limit of a subsequence $K(j)$ 
to infinity in \eqref{EQ:EnergyEQ} to follow
\allowdisplaybreaks\begin{align}
&\lint_0^t\epsilon^{1-\kappa}\,ds
\lint_{\mathbb{R}^{N+1}}\left|\frac{\partial v_{\epsilon}}{\partial t}\right|^2\, dx
\,+\,\lint_{\mathbb{R}^{N+1}}\,d\mu_t^{\epsilon}
\,+\,\lint_0^t\frac{\dot\kappa \, ds}{K(\epsilon)}\lint_{\mathbb{R}^{N+1}}\, d\xi_s^\epsilon
\notag\\                
& \,+\,\lint_0^t\frac{\ddot\kappa\,+\,\dot\kappa^2|\log\epsilon|}{K(\epsilon)\dot\kappa}
\epsilon^{1-\kappa}\, ds\lint_{\mathbb{R}^{N+1}} G(v_{\epsilon})\, dx
\;=\;\lint_{\mathbb{R}^{N+1}}\, d\mu_0^{\epsilon}
\label{EQ:EnergyEQ-SpaceTime}\end{align}
for any time \begin{math}t\,\in\, [0,t_\epsilon],\end{math} thereby completing the proof.

%
%
\setcounter{chapternumber}{4}\setcounter{equation}{0}
\renewcommand{\theequation}%
             {\thechapternumber.\arabic{equation}}
\section{\enspace Asymptotic behavior}
%
%
\par
We state an asymptotic behavior of $v_\epsilon$; 
This is necessary to prove the clearing-out lemma in a construction of Brakke's motion.
See Ilmanen \cite[\S 6.4]{ilmanen}.
\begin{Thm}{}\label{THM:Comparison}
Set \begin{math}r(x,t)\;=\;2\sqrt{T-t}\,-\,|x|\end{math} $(0\le t<T)$
for some positive number $T$\rm{;}
Let us assume that \begin{math}v_\epsilon (x,0)\;\ge\;q(r(x,0)/\epsilon).\end{math}
Then for any positive number $M$,
\allowdisplaybreaks\begin{align}
& v_\epsilon (x,t)\;\ge\;\frac{1\,-\,\epsilon^{-2M}}{1\,+\,\epsilon^{-2M}}\;\mathrm{in}\;\mathrm{any}\;
x\,\in\,\{x\in\mathbb{R}^{N+1}\,;\,r(x,t)\,>\,M\epsilon^{1-\kappa} |\log\epsilon|\}\notag\\
\intertext{and}
& v_\epsilon (x,t)\;\le\;-\frac{1\,-\,\epsilon^{-2M}}{1\,+\,\epsilon^{-2M}}\;\mathrm{in}\;\mathrm{any}\;
x\,\in\,\{x\in\mathbb{R}^{N+1}\,;\,r(x,t)\,<\,-M\epsilon^{1-\kappa}|\log\epsilon|\}\notag\end{align}
respectively holds
for any \begin{math}t\,\in\, (0,t_2)\end{math} as long as 
\allowdisplaybreaks\begin{align}
& v_\epsilon (x,0)\;>\;0 \;\mathrm{in}\;\mathrm{any}\;
x\,\in\,\{x\in\mathbb{R}^{N+1}\,;\, r(x,0)\,>\,0\}\notag\\
\intertext{and}
& v_\epsilon (x,0)\;<\;0 \;\mathrm{in}\;\mathrm{any}\;
x\,\in\,\{x\in\mathbb{R}^{N+1}\,;\, r(x,0)\,<\, 0\}\notag\end{align}
where $t_2$ is a positive number with
\begin{math}T\;=\; t_2\,+\,M^2\epsilon^{2(1-\kappa(t_2))}|\log\epsilon|^2/4.\end{math}\end{Thm}
\vskip 9pt\begin{center}{\underbar{{Proof of Theorem \ref{THM:Comparison}}}}\end{center}\vskip 6pt
\par We consider the function \begin{math}w_\epsilon (x,t)\;=\;q(r(x,t)/\epsilon^{1-\kappa})\end{math}
with \begin{math}q\,=\,\tanh.\end{math}
\par A direct computation reads
\allowdisplaybreaks\begin{align}
& \frac{\partial w_\epsilon}{\partial t}\,-\,\triangle w_\epsilon
\,+\,\frac 2{\epsilon^{2(1-\kappa)}} f(w_\epsilon)
\,-\,\frac {\dot\kappa}{K(\epsilon)} g(w_\epsilon) \notag\\
& \;=\;\frac{\dot{q}}{\epsilon^{1-\kappa}}
\left[\frac{\partial r}{\partial t}\,-\,\triangle r
\,+\,\frac {2q}{\epsilon^{1-\kappa}}(|Dr|^2\,-\,1)\right]\;\le\;0
\quad (x\ne 0).\label{INEQ:W}\end{align}
Applying \eqref{EQ:NewDis} and \eqref{INEQ:W} to a comparison theorem, we conclude
\begin{math}v_\epsilon (x,t)\,\ge\,w_\epsilon (x,t)\end{math}
in \begin{math}\mathbb{R}^{N+1}\times [0,t_\epsilon]\end{math}
because \begin{math}v_\epsilon(x,0)\,\ge\,w_\epsilon(x,0)\end{math}
and the continuity of $r$ at $x=0$.

%
%
\setcounter{chapternumber}{3}\setcounter{equation}{0}
\renewcommand{\theequation}%
             {\thechapternumber.\arabic{equation}}
\section{\enspace Brakke's motion}
%
%
\par
We finally mention the existence of Brakke's motion without the proof: 
\begin{Thm}{\rm{(}Brakke's motion\rm{)}}\label{THM:Brakke}
Assume that $g$ satisfies (i)-(iv) of p.423 Ilmanen\cite[\S 1.4]{ilmanen}
in which we respectively substitutes $g$ and $N$ for $u_0$ and $n-1$.
By employing Theorem \ref{THM:Main} and Theorem \ref{THM:Comparison} 
instead of \S 4 and \S 6.5 in Ilmanen\cite[\S 1.4]{ilmanen}, 
we can build a sequence of Radon measure $\{\mu_t\}$ $(t\ge 0)$ in definition \ref{DEF:BrakkeMotion}.
\end{Thm}

%
%
%
\bibliographystyle{amsalpha}

\end{document}